\documentclass{amsart}
\usepackage{amsmath, amsfonts, amssymb, latexsym}
\newtheorem{theorem}{Theorem}
\newtheorem{lemma}{Lemma}

\newtheorem{notation}[lemma]{Notation}
\newtheorem{problem}[lemma]{Problem}

\newtheorem{remark}[lemma]{Remark}

\begin{document}
\title{Branching Brownian motion with ``mild'' Poissonian   obstacles}
\author{J\'{a}nos Engl\"{a}nder}
\address{Department of Statistics and Applied Probability\\
University of California, Santa Barbara\\
CA 93106-3110, USA.}
\email{englander@pstat.ucsb.edu}
\urladdr{http://www.pstat.ucsb.edu/faculty/englander}
\date{\today}

\begin{abstract}
We study a spatial branching model, where the underlying motion is
 Brownian motion and the branching is
affected by a random collection of reproduction blocking sets called
{\it mild obstacles}. We  show that the quenched local growth rate
is given by the branching rate in the `free' region . When the
underlying motion is an arbitrary diffusion process, we obtain a
dichotomy for the local growth that is independent of the Poissonian
intensity. Finally, and most importantly,  we obtain the asymptotics
(in probability) of the  quenched (when $d\le 2$) and the annealed
(arbitrary $d$) {\it global growth rates}, and identify
subexponential correction terms.

\end{abstract}
\maketitle
\tableofcontents
\section{Introduction}
\subsection{Model.}
The purpose of this note is to study a spatial
branching model with the property that branching
only takes place in a certain random region  where
it is spatially (and temporarily) constant. More
specifically we will use a natural model for the
random environment: {\it mild Poissonian
obstacles}.

Let $\omega$ be a Poisson point process (PPP) on $\mathbb R^d$
with intensity $\nu
>0$ and let $\mathbf P$ denote the corresponding law.
Furthermore, let $a,\ \beta >0$ fixed. We define
the {\it branching Brownian motion (BBM) with a
mild Poissonian obstacle}, or the `$(\nu ,\beta, a
)$-BBM' as follows. Let $K$ denote the random set
given by the $a$-neighborhood of $\omega$:
$$K=K_{\omega}:=\bigcup_{x_i\in \mathrm{supp}(\omega)} \bar B(x_i,a).$$
Then $K$ is  a \emph{mild obstacle configuration} attached to
$\omega$.  This means that given $\omega$, we define
$P^{\omega}$ as the law of  the (strictly dyadic) BBM on
$\mathbb R^d,\ d\ge 1$ with spatially dependent branching rate
$r(x,\omega):=\beta 1_{K^c_{\omega}} (x)$. An equivalent
(informal) definition is that as long as a particle is in
$K^c$, it obeys the branching rule with rate $\beta$, while in
$K$ its reproduction is completely suppressed and  it does not
branch at all. (We assume that the process starts with a
single particle at the origin.) The process $Z$ under
$P^{\omega}$ is called a {\it BBM with  mild Poissonian
obstacles}. The total mass process will be denoted by $|Z|$.
Further, $W$ will denote $d$-dimensional Brownian motion with
probabilities $\{\mathbb P_{x},\ x\in\mathbb R^d \}$. Finally
$B(z,r)$ will denote the open ball centered at $z\in\mathbb
R^d$ with radius $r$ and  $\lfloor z\rfloor$ will denote the
integer part of $z\in \mathbb R$: $z:=\max \{n\in\mathbb Z\mid
n\le z \}.$

\begin{remark}[Coupling]\label{coupling}\rm
The following {\it coupling} between an ordinary
``free'' BBM and the one with mild obstacles may be
useful to have in mind. Let us randomly ``trim''
the free historical BBM tree by looking at those
branching points that are inside $K$ and delete
completely one of the two branches emanating from
that point -- the branch deleted is being chosen at
random. This way we obtain a coupling: the randomly
``trimmed'' version of the free BBM has the same
distribution as the BBM with mild obstacles.
\end{remark}

\subsection{Motivation.}\label{motiv}
Topics concerning the {\em survival asymptotics}
for a single Brownian particle {\em among
Poissonian obstacles} became a classic subject in
the last twenty years, initiated by Sznitman,
Bolthausen and others. It was originally motivated
by the celebrated `Wiener-sausage asymptotics' of
Donsker and Varadhan in 1975 and by mathematical
physics, and the available literature today is huge
--- see the fundamental monograph \cite{SZ} and the references therein.

In \cite{E} a  model of a {\it branching process in
a random environment}  has been introduced. There,
instead of  mild  obstacles, we introduced hard
obstacles and instantaneous killing of the
branching process once any particle hits the trap
configuration $K$.

In \cite{EdH} the model was further studied and a
similar model with `individual killing rule' was
also suggested. Individual killing rule means that
only the individual particle hitting $K$ is
eliminated and so the process dies when the last
particle has been absorbed. The difference between
that model and the model we consider in the present
article is that now the killing mechanism is even
`milder' as it only inhibits temporarily the
reproduction of the individual particle but does
not eliminate the particle.

Working with a BBM in a Poissonian random environment in
\cite{E,EdH} turned out to be quite challenging and resulted
in some unexpected and intriguing  results.

An alternative view on our setting is as follows. Arguably,
the model can be viewed as a {\it catalytic} BBM as well ---
the catalytic set is then $K^c$. Catalytic spatial branching
(mostly for superprocesses though) has been the subject of
vigorous research in the last twenty years initiated by
Dawson, Fleischmann and others
--- see the survey papers \cite{K} and \cite{DF}
and references therein. In those models the individual
branching rates of particles moving in space depend on the
amount of contact between the particle (`reactant') and a
certain random medium called the catalyst. The random medium
is usually assumed to be a `thin' random set (that could even
be just one point) or another superprocess.

In more complicated models `mutually' or even `cyclically'
catalytic branching is considered. (See again \cite {DF}.)

Our model is simpler than most catalytic models as our
catalytic/blocking areas are fixed, whereas in several
catalytic models they are moving. On the other hand, while for
catalytic settings studied so far results were mostly only
qualitative we are aiming to get quite sharp {\it
quantitative} result.

For the discrete setting there is much less work available,
one example is \cite{KS}. In that paper the authors study a
branching particle system
 on $\mathbb Z^d$, whose branching is
{\it catalyzed by another autonomous particle system} on $\mathbb
Z^d$. There are two types of particles, the $A$-particles
(`catalyst') and the $B$-particles (`reactant'). They move, branch
and interact in the following way. Let $N_A(x, s)$ and $N_B(x, s)$
denote the number of $A$- [resp. $B$-]particles at $x\in \mathbb
Z^d$ and at time $s\in [0,\infty)$. (Here all $N_A(x,0)$ and
$N_B(x,0)$ with $x\in \mathbb {Z}^d$ are independent Poisson
variables with mean $\mu_A$ [resp. $\mu_B$].) Every $A$-particle
($B$-particle) performs independently a continuous-time random walk
with jump rate $D_A$ ($D_B$). In addition a $B$-particle dies at
rate $\delta$, and, when present at $x$ at time $s$, it splits into
two particles in the next $ds$ time units with probability $\beta
N_A(x, s) ds+o(ds)$. Conditionally on the system of the
$A$-particles, the jumps, deaths and splitting of the $B$-particles
are independent. The authors prove that for large $\beta$ there
exists a critical $\delta$ separating local extinction regime from
local survival regime.

A further example of the discrete setting is
\cite{AB}.

\medskip
It appears that our proposed model of a BBM with
``mild'' obstacles  has biological merit to it.
First, one immediately  has the following two {\it
biological interpretations} in mind:
\begin{enumerate}
\item {\bf Migration with unfertile areas (Population dynamics):}
Population moves in space and reproduces by binary
splitting, except at randomly located
reproduction-blocking areas.
\item {\bf Fecundity
selection (Genetics):} Reproduction and mutation
takes place. Certain randomly distributed genetic
types have low fitness:  even though they can be
obtained by mutation, they themselves are unable to
reproduce, unless mutation transforms them to
different genetic types. In genetics this
phenomenon is called `{\it fecundity selection}'.
\end{enumerate}
One question of interest is of course the (local and global)
growth rate of the population. Once one knows the global
population size, the model can be rescaled (normalized) by the
global population size, giving a population of unit mass
(somewhat similarly to the fixed size assumption in the Moran
model or many other models from theoretical biology) and then
the question becomes the {\it shape} of the population.

In the population dynamics setting this latter
question concerns whether or not there is a
preferred spatial location for the process to
populate. In the genetic setting the question is
about the existence of a certain kind of genetic
type that is preferred in the long run that lowers
the risk of  low of fecundity caused by mutating
into less fit genetics types.

Of course, the {\it genealogical structure} is a
very exciting problem to explore too. For example
it seems quite possible that for large times the
`bulk' of the population consists of descendants of
a single particle that decided to travel far enough
(resp. to mutate many times) in order to be in a
less hostile environment (resp. in high fitness
genetic type area), where she and her descendants
can reproduce freely.

For example, a related phenomenon  in marine systems (personal
communication with Chris Cosner) is when hypoxic patches form in
estuaries because of stratification of the water. The patches affect
different organisms in different ways but are detrimental to some of
them. They appear and disappear in an effectively stochastic way.
This is an actual system that has some features that correspond to
the type of assumptions built into our model.

It appears (personal communication with Bill Fagan)
that a very relevant existing ecological context in
which to place our model is the so-called
``source-sink theory''. The basic idea is that some
patches of habitat are good for a species (and
growth rate is positive) whereas other patches are
poor (and growth rate is zero or negative).
Individuals can move between patches randomly or
according to more detailed biological rules for
behavior.

Another kind of scenario where models such as the
proposed one would make sense is in systems that
are subject to periodic local disturbances
(personal communication with Chris Cosner). Those
would include forests where trees sometimes fall
creating gaps (which have various effects on
different species but may harm some) or areas of
grass or brush which are subject to occasional
fires. Again, the effects may be mixed, but the
burned areas can be expected to less suitable
habitats for at least some organisms.

\medskip
\noindent Let us now return to our mathematical
model and consider the following natural questions:
\begin{enumerate}
\item What can we say about the growth of the total
population size?
\item What are the large deviations? (E.g., what is the probability of
producing an atypically small population.)
\item What can we say about the {\it local} population growth?
\end{enumerate}
Of course, these questions make sense both in the annealed and the
quenched sense.

As far as the first question is concerned, recall that the total
population of an ordinary (free) BBM grows in expectation (and
almost surely) as $e^{\beta t}$. To see this note that in fact,
for ordinary BBM, the spatial
component  plays no role, and hence the total mass is just a
$\beta$-rate pure birth process $X$. As is well known, the limit
$N:=\lim_{t\to\infty} e^{-\beta t}X_t$ exists a.s and in mean,
and $P(0<N<\infty)=1$.

Turning to BBM with the Poissonian reproduction
blocking mechanism, how much will the absence of
branching in $K$ slow the global reproduction down?
Will it actually change the exponent $\beta$? (We
will see that even though the global reproduction
does slow down, the slowdown is captured by a
sub-exponential factor, being different for the
quenched and the annealed case.)

Consider now  the second question. Here is an
argument to show the non-triviality of the problem
and to give a motivation. Let us ask the simplest
question: what is the probability that there is no
branching at all up to time $t>0$? In order to
avoid branching the first particle has to `resist'
the branching rate $\beta$ inside $K^c$. Therefore
this question is quite similar to the survival
asymptotics for a single Brownian motion among
`soft obstacles'
--- but of course in order to prevent branching the particle seeks
for large  islands covered by $K$ rather then the
usual `clearings'. In other words, the particle now
prefers to avoid the ``Swiss cheese'' $K^c$ instead
of $K$. The second listed problem above is a
possible generalization of this (modified) soft
obstacle problem for a single particle. In fact,
the presence of branching seems to bring genuinely
new type of challenges into the analysis.

Finally, the third question will be shown to be related to the
recent paper \cite{EK} that treats branching diffusions on
Euclidean domains.

\subsection{Expected global growth and dichotomy for local growth}
Concerning the expected global growth rate we have
the following result.
\begin{theorem}[Expected global growth rate]\label{exp.thm} On a set of full $\mathbf {P}$-measure,
\begin{equation}\label{quenchedasymp}
E^{\omega} |Z_t| =\exp\left[\beta t-c(d,\nu)\frac{t}{(\log
t)^{2/d}}(1+o(1))\right],\ \mathrm{as}\ t\to\infty
\end{equation}
(quenched asymptotics), and
\begin{equation}\label{annealedasymp}
\left({\mathbf E} \otimes E^{\omega}\right) |Z_t| = \exp[\beta
t-\tilde c(d,\nu)t^{d/(d+2)}(1+o(1))],\ \mathrm{as}\ t\to\infty
\end{equation} (annealed
asymptotics), where
\begin{eqnarray*}
&&c(d,\nu):=\lambda_d\left(\frac{d}{\nu\omega_d}\right)^{-2/d},\\
&&\tilde
c(d,\nu):=(\nu\omega_d)^{2/(d+2)}\left(\frac{d+2}{2}\right)\left(\frac{2\lambda_d}{d}\right)
^{d/(d+2)},
\end{eqnarray*}
and  $\omega_d$ is the volume of the $d$-dimensional  unit ball,
while $\lambda_d$ is the principal Dirichlet eigenvalue of
$-\frac1{2}\Delta$ on it.
\end{theorem}

\begin{remark} [Interpretation of Theorem 1]
\rm Consider the annealed case first. Recall the
coupling suggested in Remark \ref{coupling}. Having
that coupling in mind, let us pretend for a moment
that we are talking about an ordinary BBM. Then at
time $t$ one has $e^{\beta t}$ particles with
probability tending to one as $t\to\infty$. For $t$
fixed take a ball $B=B(0,R)$ (here $R=R(t)$) and
let $K$ be so that $B\subset K^c$ (such a ball left
empty by $K$ is called a {\it clearing}). Consider
the expected number of particles that are confined
to $B$ up to time $t$. These particles do not feel
the blocking effect of $K$, while the other
particles may have not been born due to it.

Optimize  $R(t)$ with respect to the cost of having
such a clearing and the probability of confining a
single Brownian motion to it. This is precisely the
same optimization as for the classical
Wiener-sausage. Hence one gets the expectation in
the theorem as a lower estimate.

One suspects that the main contribution in the
expectation in (\ref{annealedasymp}) is coming from
the expectation on the event of having a clearing
with optimal radius $R(t)$. In other words,
denoting by $p_t$ the probability that a single
Brownian particle stays in the $R(t)$-ball up to
time $t$, one argues heuristically that $p_t
e^{\beta t}$ particles will stay inside the
clearing up to time $t$ `for free' (i.e. with
probability tending to one as $t\to\infty$).

The intuitive reasoning is as follows. If we had {\it independent}
particles instead of BBM, then, by a `Law of Large Numbers type
argument' (using Chebysev inequality and the fact that
$\lim_{t\to\infty}p_t e^{\beta}=\infty$), roughly $p_t e^{\beta
t}$ particles out of the total $e^{\beta t}$  would stay in the
$R(t)$-ball up to time $t$ with probability tending to $1$ as
$t\uparrow \infty$.  One suspects then that the lower estimate
remains valid for the branching system too, because the particles
are ``not too much correlated'' . This kind of argument (in the
quenched case though) will be made precise in the proof of our
main theorem by estimating certain covariances.

The interpretation of the formula for the quenched case is
different in the sense that large clearings (far away) are
automatically (that is, $\mathbf{P}$-a.s.) present; but it is
similar in the sense that part of our job will be to show that
inside such a clearing a large population is going to flourish
(see the proof of Theorem \ref{mainthm} for more on this). $\hfill
\Diamond$
\end{remark}

Concerning local population size we have the following (quenched)
result.
\begin{theorem}[Quenched exponential growth]\label{qeg.blo}
The following holds on a set of full $\mathbf {P}$-measure: For
any $\epsilon>0$ and any bounded open set $\emptyset\neq B\subset
\mathbb {R}^d$,
\begin{equation*}
P^{\omega}_{\mu }\left( \limsup_{t\uparrow \infty
}e^{-(\beta-\epsilon ) t}Z_{t}(B)=\infty \right) >0 \mbox{ and
}P^{\omega}_{\mu }\left( \limsup_{t\uparrow \infty }e^{-\beta
t}Z_{t}(B)<\infty \right) =1.
\end{equation*}
\end{theorem}

\begin{problem} \rm What can one say about the {\em distribution} of the
global and the local population size? Our theorems make it
plausible  that  the global population size is a random multiple
of  $\exp\left[\beta t-c(d,\nu)\frac{t}{(\log
t)^{2/d}}(1+o(1))\right]$ as $t\uparrow \infty$ for almost every
environment. This problem will be addressed in Section \ref{LLN}.
We will  prove an appropriate formulation of the statement when
the limit is meant {\it in probability} and $d\le 2$.
\end{problem}
We now show how Theorem \ref{qeg.blo} can be generalized for the
case when the underlying motion is a diffusion. Let $\mathbf P$ be
as before but replace the Brownian motion by an $L$-diffusion $X$
on $\mathbb R^d$, where $L$ is a second order elliptic operator of
the form $$L=\sum_{i,j=1}^na_{i,j}(x)\frac{\partial^2}{\partial
x_i
\partial x_j}+\sum_{i=1}^nb_i(x)\frac\partial{\partial x_i}$$
with $a_{i,j}, b_i\in C^\alpha(R^n)$, $\alpha\in [0,1)$, and
the symmetric matrix $a_{i,j}(x)$ is positive definite for all
$x\in \mathbb R^d$. The branching $L$-diffusion with the
Poissonian obstacles can be defined analogously to the case of
BM. To present the result, we need an additional concept. Let
$$ \lambda _{c}(L)=\lambda _{c}(L ,\mathbb R^d):=\inf
\{\lambda \in \mathbb{R}\ :\ \exists u>0\ \text{satisfying}\
(L -\lambda )u=0\ \text{in}\ \mathbb R^d\} $$ denote the
\textit{generalized principal eigenvalue }for $L$ on $\mathbb
R^d$. In fact $\lambda _{c}\le 0$, because $L1=0$.
--- see \cite{P}, section 4.4.

The following theorem shows that the local behavior of the process
exhibits a dichotomy. The crossover is given in terms of the local
branching rate
 $\beta$ and the term $\lambda_c (L)$: local extinction occurs when the branching
 rate inside the `free region' $K^c$ is not sufficiently large to compensate
 the transience of the underlying $L$-diffusion; if it is strong
 enough, then local mass grows exponentially.
Note an interesting feature of the result: {\em the
intensity $\nu$ of the obstacles plays no role }.
\begin{theorem}[Quenched exponential growth/local extinction]\label{qegle.blo}
Given the environment $\omega$, denote by $P^{\omega}$  the law of
the branching $L$-diffusion.
\begin{itemize}
\item[(i)] Let $\beta>-\lambda_c (L)$ and let $\nu>0$ be arbitrary.
Then the following holds on a set of full $\mathbf {P}$-measure:
For any $\epsilon>0$ and any bounded open set $\emptyset\neq
B\subset \mathbb {R}^d$,
\begin{equation*}
P^{\omega}\left( \limsup_{t\uparrow \infty
}e^{(-\beta-\lambda_c (L)+\epsilon )
t}Z_{t}(B)=\infty \right) >0.
\end{equation*}
and
\begin{equation*}
P^{\omega}\left( \limsup_{t\uparrow \infty
}e^{(-\beta-\lambda_c (L)) t}Z_{t}(B)<\infty
\right) =1.
\end{equation*}
\item[(ii)] Let $\beta\le - \lambda_c (L)$ and let $\nu>0$ be
arbitrary. Then the following holds on a set of full $\mathbf
P$-measure: For any  bounded open set $ B\subset \mathbb {R}^d$
there exists a $P^{\omega}$-a.s. finite random time $t_0=t_0(B)$
such that $X_t(B)=0$  for all  $t\ge t_0$,  (local extinction).
\end{itemize}
\end{theorem}
\subsection{Quenched asymptotics of  global growth}\label{LLN}
In this section we assume that $d\le 2$ and investigate the behavior
of the (quenched) global growth rate. Define the {\it average growth
rate} by
$$r_t=r_t(\omega):=\frac{\log|Z_t(\omega)|}{t}.$$
Replace now $|Z_t(\omega)|$ by its expectation
$\bar{Z}_t:=E^{\omega}|Z_t(\omega)|$, and define
$$\hat{r}_t=\hat {r}_t(\omega):=\frac{\log \bar{Z}_t}{t}.$$
Recall
from Theorem \ref{quenchedasymp}, that on a set of
full $\mathbf {P}$-measure,
\begin{equation}
\lim_{t\to\infty}(\log t)^{2/d}(\hat
{r}_t-\beta)=-c(d,\nu ).
\end{equation}
We are going to show that an analogous statement
holds for $r_t$ itself.
\begin{theorem} \label{mainthm} Let $d\le 2$. On a set of full $\mathbf
{P}$-measure,
\begin{equation}
\lim_{t\to\infty}(\log t)^{2/d}(r_t-\beta)=-c(d,\nu
)\ \ \ \ \ \  \mathrm{in\ P^{\omega}-probability}.
\end{equation}
\end{theorem}
(That is, loosely speaking, $r_t\approx \beta-c(d,\nu)(\log
t)^{-2/d}.)$
\subsection{Outline}
The rest of this article is organized as follows.  Section 2
discusses some further problems. Section 3 presents some
preparatory results, while section 4 contains the proofs of the
theorems. Section 5 presents some additional problems that refer
to certain parts of the proofs. In section 6 we discuss yet
another result. Finally the appendix presents a proof of a
statement that is needed in the proofs and that we consider to be
of independent interest.
\section{Further problems}
In this section we suggest some further problems and
directions for research.
\subsection{More general branching} It should also
be investigated, what happens when dyadic branching
is replaced by a {\it general} one. In a more
sophisticated population model, particles can also
die --- then the obstacles {\it do not necessarily
reduce the population size} as they sometimes
prevent death.

For example {\it critical branching} requires an
approach very different from the supercritical one,
since taking expectation now does not  provide a
clue: $E^{\omega}|Z_t(\omega)|=1,$ $\forall t>0,\
\forall\omega\in \Omega$.

It is not even clear  that
$$P^{\omega}(\mathrm{extinction})=1\
\forall\omega\in \Omega .$$ And if it is true, then what can
we say  about the asymptotic behavior (as $t\to\infty$) of
$P(\tau>t)$? (Here $\tau$ denotes extinction time.)

\subsection{Superprocesses with mild obstacles}

A further goal is to generalize the setting by
defining {\it superprocesses with mild obstacles}
analogously to the BBM with mild obstacles. Recall
the concept of an
$(L,\beta,\alpha)$--superdiffusion: Let
$\mathcal{M}_{\mathrm{f}}$ denote the set of finite
measures $\mu$ on $\mathbb R^d $ and let
$\alpha,\beta$ denote functions in the H\"{o}lder
space $C^{\gamma}$ satisfying $ \alpha>0\
\text{and}\ \sup_{x\in \mathbb R^d}\beta(x)<\infty.
$
\begin{notation}\label{notat}
\emph{Let }$\left(  X,\mathbf{P}_{\mu\,},\,\mu\in
\mathcal{M}_{\mathrm{f}}\right)  $ \emph{denote the
}$(L,\beta,\alpha)$--superdiffusion\emph{. That is,
}$X$ \emph{is the unique}
$\mathcal{M}_{\mathrm{f}}$\emph{--valued
(time-homogeneous) continuous Markov process which
satisfies, for any bounded continuous }$g:\mathbb
R^d \mapsto\mathbb{R}_{+\,}$,
\begin{equation}
\mathbf{E}_{\mu}\exp\left\langle
X_{t\,},-g\right\rangle \;=\;\exp\left(
-\int_{\mathbb R^d}\mu(\mathrm{d}x)\;u(t,x)\right)
,
\label{Laplace.functional}%
\end{equation}
\emph{where} $u$ \emph{is the minimal non-negative solution to}%
\begin{equation}
\left.
\begin{array}
[c]{c}%
\dfrac{\partial}{\partial t}u\,=\,Lu+\beta u-\alpha u^{2}\quad\text{\emph{on}%
}\;\;\mathbb R^d\times(0,\infty),\vspace{8pt}\\
\lim\limits_{t\rightarrow0+}u(t,\,\cdot\,)\,=\,g(\cdot)
\end{array}
\;\right\}  \label{cum.equ}%
\end{equation}
\rm (Here $\left\langle \nu,f\right\rangle $
denotes the integral $\int_{\mathbb
R^d}\nu(\mathrm{d} x)\,f(x).$)
\end{notation}

The definition of the superprocess with mild
obstacles is straightforward:  the parameter
$\alpha$ vanishes on the set $K$ and elsewhere it
is positive (for example, $\alpha:=\alpha_0
1_{K^c}$ where $\alpha_0$ is a positive constant).

Similarly, one can consider the case when instead of $\alpha$, the
`mass creation term' $\beta$ is random, e.g. $\beta:=\beta_0
1_{K^c},\ \beta_0>0$. Denote now by $P^{\omega}$ the law of this
latter superprocess for given environment. We suspect that the
superprocess with mild obstacles behaves similarly to the discrete
branching process with mild obstacles when $\lambda_c(L+\beta)>0$
and $P^{\omega}(\cdot)$ is replaced by $P^{\omega}(\cdot\mid X\
\mathrm{survives} )$.


\section{Preparations} In this section we present some
preparatory lemmas for the proofs.

\begin{lemma}[Expectation given by Brownian functional]
Fix $\omega$. Then
\begin{equation}\label{FK} {E}^{\omega} |Z_t| = \mathbb
E \exp \left[\int_0^{t}\beta
1_{K^c}(W_s)\right]\mathrm{d}s\end{equation}
\end{lemma}
{\bf Proof.} It is well known (`first moment
formula' of spatial branching processes) that
${E}_x^{\omega} |Z_t|=(T_t 1)(x)$, where
$u(x,t):=(T_t 1)(x)$ is the minimal solution of the
parabolic problem:
\begin{eqnarray}\label{min.par.}
\frac{\partial u}{\partial
t}&=&\left(\frac{1}{2}\Delta +\beta 1_K^c\right) u\
\mathrm{on}\ \mathbb R^d\times (0,\infty ),\nonumber\\
u(\cdot,0)&=&1,\\
u&\ge& 0.\nonumber
\end{eqnarray}
This is equivalent (by the Feynman-Kac formula)
with (\ref{FK}). \hfill $\blacksquare$

We will also need the following result :
\begin{lemma}
[\cite{EK}, Theorem 3]\label{EK}

Let  $0\le V$, $0\not\equiv V$ be bounded from above and smooth.
Let $P$ denote the law of the $(L,V)$-branching diffusion $X$ and
let $\lambda_c:=\lambda_c(L+V).$
\begin{itemize}

\item[(i)]  Under $P$ the process $X$ exhibits  local extinction if and only
if  $\lambda _{c}\leq 0$.

\item[(ii)]  When $\lambda _{c}>0$, for any $\lambda <\lambda _{c}$
and  any open $\emptyset\neq B\subset \subset D$,
\begin{equation*}
P_{\mu }\left( \limsup_{t\uparrow \infty }e^{-\lambda
t}X_{t}(B)=\infty \right) >0 \mbox{ and }P_{\mu }\left(
\limsup_{t\uparrow \infty }e^{-\lambda _{c}t}X_{t}(B)<\infty
\right) =1.
\end{equation*}
\end{itemize}
\end{lemma}

\section{Proofs of Theorems 1-\ref{mainthm}.}
\subsection{Proof of Theorem 1.}
 We can rewrite
the equation (\ref{FK}) as
$${E}^{\omega} |Z_t| = e^{\beta t} \mathbb E \exp \left[-\int_0^{t}\beta
1_{K}(W_s)\right]\mathrm{d}s.
$$
The expectation on the righthand side is precisely the survival
probability among `soft obstacles', except that we do not sum the
shape functions on the overlapping balls. However this does not
make any difference with regard to the asymptotics (see \cite{SZ},
Remark 4.2.2.). The statements thus follow from  the well known
Donsker-Varadhan type `Wiener asymptotics' for soft obstacles (
\cite{SZ} Theorems 4.5.1. and 4.5.3.; see also see \cite{DV}).
$\hfill \blacksquare$

\subsection{Proof of Theorem \ref{qeg.blo}.} Since in this case
$L=\Delta/2$ and since $\lambda_c(\Delta/2,\mathbb R^d)=0$, the
statement immediately follows from Theorem \ref{qegle.blo}.
$\hfill \blacksquare$

\subsection{Proof of Theorem \ref{qegle.blo}.}  Let  $g$ and $G$
be two functions from $C^{\gamma}(\mathbb R^d)$ (here $\gamma\in(0,1]$)
with $$1_{B(0,a/2)}\le g\le
1_{B(0,a)}\le G\le 1,$$ and let
$$V_1(\cdot):=\beta\sum_{\omega}g(\cdot+\omega),\ V_2(\cdot):=\beta\sum_{\omega}G(\cdot+\omega).$$
Then \begin{equation} V_1\le\beta 1_{K^c}\le
V_2.\label{sandwich}
\end{equation}
Consider the operators $L+V_i$ ($i=1$ or $i=2$) on
$\mathbb R^d$ and let $\lambda^i_c=\lambda^i_c
(\omega )$ denote their generalized principal
eigenvalues. Note that $V_i\in C^{\gamma}(\mathbb R^d)$
and thus
standard theory is applicable. In particular, since $V_1\le V_2$, one has
$\lambda^1_c\le\lambda^2_c$ (see Chapter 4 in
\cite{P}). As a first step, we show that in fact
\begin{equation}\label{princeig}
\lambda^1_c=\lambda^2_c= \lambda_c(L)+\beta
\mathrm{\ for\ a.e.\ }\omega.
\end{equation}
By obvious comparison with the operator $L+\beta $
on $\mathbb R^d$, $\lambda^2_c \le
\lambda_c(L)+\beta$ for every $\omega$. On the
other hand, one gets a lower estimate on
$\lambda^1_c$ as follows. Fix $R>0$. By the
homogeneity of the Poisson point process, for
almost every environment,  supp$V_1$ contains a
clearing of radius $R$. Hence, by comparison,
$\lambda^1_c\ge \lambda^{(R)}$, where
$\lambda^{(R)}$ is the principal Dirichlet
eigenvalue of $L+\beta$ on a ball of radius $R$.
Since $R$ can be chosen arbitrarily large and since
$\lim_{R\uparrow
\infty}\lambda^{(R)}=\lambda_c(L)+\beta$, we
conclude that $\lambda^1_c\ge \lambda_c(L)+\beta$
for almost every environment.

From (\ref{sandwich}) and (\ref{princeig}) it
follows that
\begin{equation}\label{orig.princeig}
\lambda(L+\beta 1_{K^c})= \lambda_c(L)+\beta
\mathrm{\ for\ a.e.\ }\omega.
\end{equation}
The statements of the theorem  follow from (\ref{orig.princeig})
along with Lemma \ref{EK}. $\hfill \blacksquare$

\subsection{Proof of Theorem \ref{mainthm}.}

We give an upper and a lower estimate separately.

\medskip
\noindent\underline{Upper estimate:} Let
$\epsilon>0.$ Using  the Markov and
Jensen-inequalities along with the expectation
formula (\ref{quenchedasymp}), we have that on a
set of full $\mathbf {P}$-measure:
\begin{eqnarray*} && P^{\omega}\left((\log
t)^{2/d}(r_t-\beta)+c(d,\nu )>\epsilon\right)\le
\epsilon^{-1}[(\log t)^{2/d}(E^{\omega}
r_t-\beta)+c(d,\nu )]\\
&&\le\epsilon^{-1}\left[(\log
t)^{2/d}\left(\frac{\log
E^{\omega}|Z_t|}{t}-\beta\right)+c(d,\nu
)\right]=o(1)\ \mathrm{as}\ t\to\infty.
\end{eqnarray*}

\noindent\underline{Lower estimate:} We give a
``bootstrap argument'': we first prove a weaker
estimate which we will later use to prove the
stronger result. The weaker result is as follows.
Let $0<\delta<\beta$. Then on a set of full
$\mathbf{P}$-measure
\begin{equation}\label{delta}
\lim_{t\to\infty}P^{\omega}(|Z_t|\ge e^{\delta
t})=1.
\end{equation}
To prove (\ref{delta}), recall (\ref{princeig}) and that
$\lambda_c(\Delta)=0$. Using these, take $R>0$ large enough so that
$\lambda_c=\lambda_c(B_R(0))$, the principal eigenvalue of the
linearized operator $\frac{1}{2}\Delta +\beta 1_{K^c}$  satisfies
$$\lambda_c(B_R(0))>\delta.$$
Let $\hat Z^R$ be the process obtained from $Z$ by introducing
killing at $\partial B_R(0)$. Then
$$\lim_{t\to\infty}P(|Z_t|<e^{\delta t})\le \lim_{t\to\infty}
P(|\hat Z^R_t|<e^{\delta t}).$$ Let $0\le \phi=\phi^R$ be the
Dirichlet eigenfunction corresponding to $\lambda_c$ on $B_R(0)$,
and normalize it by $\sup_{x\in B_R(0)}\phi (x)=1$. Then we can
continue with
$$\le \lim_{t\to\infty}
P(\langle \hat Z^R_t, \phi\rangle<e^{\delta t}),$$ where $\langle
\hat Z^R_t, \phi\rangle:=\sum_{i}\phi (\hat Z^{R,i}_t)$ and $\{\hat
Z^{R,i}_t \}$ is the `$i$th particle' in $\hat Z^R_t$. Notice that
$M_t=M^{\phi}_t:=e^{-\lambda_c t}\langle \hat Z^R_t, \phi\rangle$ is
a non-negative martingale (see p. 84 in \cite{EK}), and define
$$N:=\lim_{t\to\infty}M_t.$$
Since $\lambda_c(B_R(0))>\delta,$ the estimate is then continued as
\begin{align}\label{N}=\lim_{t\to\infty}
P(M_t<e^{(\delta-\lambda_c) t}\mid N=0)\ P^{(R)}(N=0)\le
P^{(R)}(N=0).\nonumber
\end{align}
We have that
$$\lim_{t\to\infty}P(|Z_t|<e^{\delta t})\le
P^{(R)}(N=0)$$ holds for {\it all} $R$ large enough. Therefore, in
order to prove (\ref{delta}), we have to show that
\begin{equation}
\lim_{R\to\infty} P^{(R)}(N>0 )=1.
\end{equation}
Consider now the elliptic boundary value problem
\begin{equation}
\begin{array}{l}
\displaystyle \frac{1}{2}\Delta
u+\beta 1_{K^c}(u-u^2)=0\ \mathrm{in}\ B_R(0), \\
\lim_{x\rightarrow \partial B_R(0)}u(x)=0, \\
u>0\ \mathrm{in}\ B_R(0).\label{bvp}
\end{array}
\end{equation}
The existence of a solution follows from the fact that $\lambda_c>0$
by an analytical argument given in \cite{Pconst}. (Uniqueness
follows by the semilinear maximum principle.)

The   argument  below  gives a {\it probabilistic construction} for
the solution. Namely, we show that $w_R(x):=P^{(R)}_x(N>0)$ solves
(\ref{bvp}). To see this, let $v=v_R:=1-w_R$. Let us fix an
arbitrary time $t>0$. Using the Markov and branching properties of
$Z$ at time $t$, it is straightforward to show that
$$P(N=0\mid\mathcal{F}_t)=\prod_i P_{\hat Z^{R,i}_t}(N=0).$$
Since the left hand side of this equation defines a martingale in
$t$, so does the right hand side. That is
$$\widetilde{M}_t:=\prod_i v(\hat Z^{R,i}_t)$$ defines a martingale.
From this, it follows by Theorem 17 of \cite{EKpreprint} that  $v$
solves the  equation obtained from the first equation of (\ref{bvp})
by switching $u-u^2$ to $u^2-u$. Consequently,
$w_R(x):=P^{(R)}_x(N>0)$ solves the first equation of (\ref{bvp})
itself. That $w_R$ solves the second equation, follows easily from
the continuity of Brownian motion. Finally its positivity (the third
equation of (\ref{bvp})) follows again from the fact that
$\lambda_c>0$ (see Lemma 6 in \cite{EK}).

By the semilinear elliptic maximum principle (Proposition 7.1 in
\cite{EP99}; see also \cite{Pconst}), $w_R(\cdot)$ is monotone
increasing in $R$. Using standard arguments, one can show that
$0<w:=\lim_{R\to\infty} w_R$ solves (\ref{bvp}) too (see the proof
of Theorem 1 in \cite{Pconst}).

Applying the strong maximum principle to $v:=1-w$, it follows that
$w$ is either one everywhere or less than 1 everywhere. Now suppose
that $0<w<1$. Then
$$\frac{1}{2}\Delta w=\beta 1_{K^c}(w^2-w)\lneqq 0.$$
This contradicts the recurrence of the Brownian motion in one and
two dimensions (see Chapter 4 in \cite{P}) . This contradiction
proves that in fact $w=1$ and consequently it proves  (\ref{delta}).

We have now completed the first part of our ``bootstrap'' proof.

\medskip
Let us return to the proof of the lower estimate in Theorem
\ref{mainthm}. Let $\epsilon>0$. We have to show that on a set of
full $\mathbf {P}$-measure,
\begin{equation}\label{need}
\lim_{t\to\infty}P^{\omega}\left((\log
t)^{2/d}(r_t-\beta)+c(d,\nu )<-\epsilon\right)=0.
\end{equation}
To achieve this, we will define a particular function $p_t$
(the definition is given in (\ref{pt})) satisfying that as
$t\to\infty$,
\begin{equation}
p_t=\exp\left[-c(d,\nu)\frac{t}{(\log t)^{2/d}}+
o\left(\frac{t}{(\log t)^{2/d}}\right)\right].
\end{equation}
Using this function we are going to show a statement implying
(\ref{need}), namely, that for all $\epsilon>0$ there is a set of
full $\mathbf {P}$-measure, where
\begin{equation}\label{stronger}
\lim_{t\to\infty}P^{\omega}\left(\log|Z_t|<\beta t+\log p_t-\epsilon
t(\log t)^{-2/d}\right)=0.
\end{equation}
Let us first give an { outline of the {\it strategy of our proof}.
A key step will be introducing three different time scales,
$\ell(t)$, $m(t)$ and $t$ where $\ell(t)=o(m(t))$ and $m(t)=o(t)$
as $t\to\infty$. For the first, shortest time interval, we will
use that there are ``many'' particles produced and they are not
moving ``too far away'', for the second (of length $m(t)-\ell
(t)$) we will use that one particle moves into a clearing of a
certain size at a certain distance, and in the third one (of
length $t-m(t)$) we will use that there is a  branching tree
emanating from that particle so that a certain proportion of
particles of that tree stay in the clearing with probability
tending to one.

To carry out this program, first recall the
following fact (for example this can be found in
the proof of Theorem 4.5.1 in \cite{SZ}): Let
$$R_0=R_0(d,\nu):=\sqrt{\frac{\lambda_d}{c(d,\nu)}}=
\left(\frac{d}{\nu\omega_d}\right)^{1/d},$$ (recall that $\lambda_d$
is the principal Dirichlet eigenvalue of $-\frac{1}{2}\Delta$ on the
$d$-dimensional unit ball, and $\omega_d$ is the volume of that
ball) and let
\begin{equation}\label{rho}
\rho=\rho(\ell):=R_0(\log \ell )^{1/d}-(\log\log\ell)^2,\ t\ge
0.
\end{equation} Then,
\begin{equation}\label{clearing}
\mathbf{P}(\exists\ \ell_0(\omega)>0\ \ \mathrm{such\ that}\ \
\forall \ell>\ell_0(\omega)\ \exists\ \mathrm{clearing}\
B(x_0,\rho)\ \mathrm{with}\ |x_0|\le\ell)=1.
\end{equation}
  Next, let $\ell$ and $m$ be two functions $\mathbb
R_+\rightarrow \mathbb R_+$ satisfying the following:
\begin{enumerate}
\item[(i)] $$\lim_{t\to\infty} \ell(t)=\infty,$$

\item [(ii)] $$\lim_{t\to\infty}\frac{\log t}{\log \ell(t)}=1,$$

\item [(iii)] $$\ell(t)=o(m(t)),\ \mathrm{as}\  t\to\infty$$

\item[(iv)] $$m(t)=o(\ell^2(t)),\ \mathrm{as}\  t\to\infty$$

\item[(v)] $$m(t)=o(t(\log t)^{-2/d}),\ \mathrm{as}\ \ t\to\infty.$$

\end{enumerate}
Note that $(i)-(v)$ are in fact not independent, because $(iv)$
follows from $(ii)$ and $(v)$. For example the following choices of
$\ell$ and $m$ satisfy $(i)-(v)$: let $\ell(t)$ and $m(t)$ be
arbitrarily defined for $t\in [0,e]$, and
$$\ell(t):=t^{1-1/(\log\log t)},\ m(t):=t^{1-1/(2\log\log t)},
\ \mathrm{for}\ t\ge t_0> e.$$ Fix $\delta\in (0,\beta)$ and define
$$I(t):=\lfloor\exp(\delta \ell(t))\rfloor.$$
Let $A_t$ denote the following event:
$$A_t:=\{|Z_{\ell(t)}|\ge I(t) )\}.$$
By (\ref{delta}) we know that on a set of full
$\mathbf{P}$-measure,
\begin{equation}\label{manyatell(t)}
\lim_{t\to\infty}P^{\omega}(A_t)=1.
\end{equation}
By (\ref{manyatell(t)}), for $t$ fixed  we can work on
$A_t\subset\Omega$ and consider $I(t)$ many particles at time
$\ell(t)$.

As a next step, we need some control on their spatial position. To
achieve this, use Remark \ref{coupling} to compare BBM's with and
without obstacles, and then, use the following result taken from
Proposition 2.3 in \cite{EdH} (the stronger, a.s. result is proved
in \cite{Ky}).

Let (only for the following statement) denote $Z$ the BBM {\it
without} obstacles starting at the origin with a single
particle. Let $R(t) = \cup_{s\in[0,t]}\, \mathrm{supp}(Z(s))$
denote the range of $Z$ up to time $t$. Let
\begin{equation}\label{radialspeed}
M(t) = \inf \{r>0: R(t) \subseteq B(0,r)\}\ \mathrm{for}\ d
\geq 1,
\end{equation}
be the radius of the minimal ball containing $R(t)$. Then
$M(t)/t$ converges to $\sqrt{2\beta}$ in probability as
$t\to\infty$.

Going back to the set of $I(t)$ many particles at time
$\ell(t)$, (\ref{radialspeed}) yields that even though they
are at different locations, still for any $\epsilon'>0$, with
$P^{\omega}$-probability tending to one, they are all inside
the $(\sqrt{2\beta}+\epsilon')\ell(t)$-ball.

Define
$$\rho(t):=R_0[\log \ell (t)]^{1/d}-[\log\log\ell(t)]^2,\ t\ge 0.$$
Recall (\ref{clearing}).  With $\mathbf{P}$-probability one
there is a clearing $B=B(x_0,\rho(t))$  such that $|x_0|\le
\ell(t)$, for all large enough $t>0$. In the sequel we may
assume without the loss of generality that  $|x_0|= \ell(t)$.
(Of course, $x_0$ depends on $t$, but this dependence is
suppressed in our notation.) By the previous paragraph, with
$P^{\omega}$-probability tending to one, the distance of $x_0$
from {\it each} of the $I(t)$ many particles is at most
$$(1+\sqrt{2\beta}+\epsilon') \ell(t).$$
Now, any such particle moves to $B(x_0,1)$ in another
$m(t)-\ell(t)$ time with probability $q_t$, where (using
$(iii)$ and $(iv)$)
$$q_t=\exp\left(-\frac{[(1+\sqrt{2\beta}+\epsilon')\ell(t)]^2}
{2[m(t)-\ell(t)]}+o\left(\frac{[(1+\sqrt{2\beta}+\epsilon')\ell(t)]^2}
{2[m(t)-\ell(t)]}\right)\right)\to 0\ \ \ \mathrm{as}\
t\to\infty.$$   Let the particle positions at time $\ell(t)$ be
$z_1,z_2,...,z_{I(t)}$ and consider the independent system of
Brownian particles $$\{W_{z_i};\ i=1,2,...,I(t)\},$$ where
$W_{z_i}(0)=z_i;\ i=1,2,...,I(t)$. In other words, $\{W_{z_i};\
i=1,2,...,I(t)\}$ just describes the evolution of the $I(t)$
particles picked at time $\ell(t)$ without respect to their
possible further descendants and (using the Markov property) by
resetting the clock at time $\ell(t)$ .

Let $C_t$ denote the following event:
$$C_t:=\{\exists i\in \{1,2,..., I(t)\},\ \exists\, 0\le s\le m(t)-\ell(t)
\ \mathrm{such\ that}\  W_{z_i}(s)\in B(x_0,1)\}.$$ By the
independence of the particles,
\begin{equation}\label{onegoesthere}
\limsup_{t\to\infty}{P}^{\omega}\left(C^c_t\mid
A_t\right)=\limsup_{t\to\infty}(1-q_t)^{
I(t)}=\limsup_{t\to\infty}\left[(1-q_t)^{1/q_t}\right]^{q_t
I(t)}.
\end{equation}
Since $(iii)$ implies that
$\frac{\ell^2(t)}{m(t)}=o(\ell(t))\  \mathrm{as}\
t\to\infty$ and since $(i)$ is assumed, one has
$$q_te^{\delta\ell(t)}=\exp\left(-\frac{[\ell(t)+(\sqrt{2\beta}+\epsilon')\ell(t)]^2}
{2[m(t)-\ell(t)]}+\delta
\ell(t)+o(\ell(t))\right)\to \infty\ \ \
\mathrm{as}\ t\to\infty.$$ In view of this,
(\ref{onegoesthere}) implies that
$\lim_{t\to\infty}{P}^{\omega}\left(C^c_t\mid
A_t\right)=0$. Using this  along with
(\ref{manyatell(t)}), it follows that on a set of
full $\mathbf{P}$-measure,
\begin{equation}\label{onewillmakeit}
\lim_{t\to\infty}{P}^{\omega}\left(C_t\right)=1.
\end{equation}
Once we know (\ref{onewillmakeit}), we proceed as follows. Recall
that $x_0$ denotes the center of $B$ and  that $\mathbb P$ denotes
the probability corresponding to a single Brownian particle $W$.
Let $\sigma^{x_0,t}_B$ denote the first exit time from $B$:
$$\sigma^{x_0,t}_B:=\inf\{s\ge 0\mid W_s\not\in B\}.$$ Abbreviate $t^*:=t-m(t)$
and define
\begin{equation}\label{pt}
p_t:=\sup_{x\in B(x_0,1)}\mathbb P_x(\sigma^{x_0,t}_B\ge t^*)=
\sup_{x\in B(0,1)}\mathbb P_x(\sigma_B\ge t^*),
\end{equation}
where $$\sigma_B:=\inf\{s\ge 0\mid W_s\not\in
B(0,\rho(t))\}.$$ It is easy to check that as $t\to\infty$,
\begin{equation}\label{ptfirstasympt}
p_t =\exp\left[-c(d,\nu)\frac{t^*}{(\log\ell(t))^{2/d}}+
o\left(\frac{t^*}{(\log\ell(t))^{2/d}}\right)\right]
;\end{equation} and using $(ii)$ and $(v)$, it follows that in
fact
\begin{equation}\label{ptasympt}
p_t=\exp\left[-c(d,\nu)\frac{t}{(\log t)^{2/d}}+
o\left(\frac{t}{(\log t)^{2/d}}\right)\right].
\end{equation}
A little later we will also need the following notation:
\begin{equation}\label{pst}
p_s^t:=\sup_{x\in B(x_0,1)}\mathbb P_x(\sigma^{x_0,t}_B\ge s)=
\sup_{x\in B(0,1)}\mathbb P_x(\sigma_B\ge s),
\end{equation}
 With this notation, $$p_t=p_{t^*}^t.$$

By slightly changing the notation, let $Z^{x}$ denote the BBM
starting with a single particle at $x\in B$; and let $Z^{x,B}$
denote the BBM starting with a single particle at $x\in B$ and
with absorbtion at $\partial B$ (and still branching at the
boundary at rate $\beta$).

Since branching does not depend on motion, $|Z^{x,B}|$ is a
non-spatial Yule's process (and of course it does not depend
on $x$)  and thus for all $x\in B$,
\begin{equation}\label{N}
\exists N:=\lim_{t\to\infty}e^{-\beta t}|Z_t^{x,B}|>0
\end{equation}
almost surely (see Theorems III.7.1-2 in \cite{AN}).

Recall Remark \ref{coupling}. By the coupling described there, it is
clear that one can in fact define the process $Z$ and the random
variable $N$ on the same probability space. Note that some particles
of $Z^{x}$ may re-enter $B$ after exiting, whereas for $Z^{x,B}$
that may not happen. Thus, by a simple coupling argument, one has
that for all $t\ge 0$, the random variable $|Z^x_t(B)|$ is
stochastically larger than $|Z^{x,B}_t(B)|$.

Recall  that our goal is to show (\ref{stronger}), and recall also
(\ref{pt}) and (\ref{ptasympt}). In fact, we will prove the
following, somewhat stronger version of (\ref{stronger}): we will
show that if the function $\gamma:\ [0,\infty)\rightarrow
[0,\infty)$ satisfies $\lim_{t\to\infty}\gamma_t=0$, then on a set
of full $\mathbf {P}$-measure,
\begin{equation}\label{evenstronger}
\lim_{t\to\infty}P^{\omega}\left(|Z_t|<\gamma_t \cdot e^{\beta
t^*}p_t\right)=0.
\end{equation}
Recalling  $t^*=t-m(t)$, and setting
$$\gamma_t:=\exp\left(m(t)-\epsilon \frac{t}{(\log t)^{2/d}}\right),
\ \mathrm{for}\ t\ge t_0> e,$$
a simple computation shows that (\ref{evenstronger}) yields
(\ref{stronger}).
 Note that this
particular $\gamma$ satisfies $\lim_{t\to\infty}\gamma_t=0$ because
of the condition $(v)$ on the function $m$.

By the comparison between $|Z^x_t(B)|$ and $|Z^{x,B}_t(B)|$
(discussed in the paragraph after (\ref{N})) along with
(\ref{onewillmakeit}) and the Markov property applied at time
$m(t)$, we have that
$$\lim_{t\to\infty}P\left(|Z_t|<\gamma_t
\cdot e^{\beta t^*}p_t\right)\le\lim_{t\to\infty}\sup_{x\in
B}P\left(|Z^{x,B}_{t^*}(B)|<\gamma_t \cdot e^{\beta
t^*}p_t\right).$$

Consider now the $J(x,t):=|Z^{x,B}_{t^*}|$ many (correlated)
Brownian paths  starting at $x\in B$ and let us denote them by
$W_1,...,W_{J(x,t)}$. Let
$$n_t^x:= \sum_{i=1}^{J(x,t)}1_{A_i},$$
where $$A_i:=\{W_i(s)\in B,\ \forall\ 0\le s\le t\}.$$  Then we have
 to show that
\begin{equation}
L:=\lim_{t\to\infty}\sup_{x\in B} P\left(n_t^x<\gamma_t \cdot
e^{\beta t^*}p_t\right)=0.
\end{equation}
Having the coupling between $Z$ and $N$ in mind, clearly, for all
$x\in B$,
\begin{eqnarray}\nonumber &L=\lim_{t\to\infty} \sup_{x\in
B}P\left(\frac{n_t^x}{N e^{\beta t^*}}<\frac{\gamma_t p_t}{N}
\right)\le\\ & \lim_{t\to\infty}\sup_{x\in B} \left[
P\left(\frac{n_t^x}{N e^{\beta t^*}}<\frac{1}{2}
p_t\right)+P\left(N\le 2\gamma_t\right)\right]. \end{eqnarray}
 Using the fact that
$\lim_{t\to\infty}\gamma_t=0$ and that $N$ is almost surely
positive, $$\lim_{t\to\infty}P\left(N\le 2\gamma_t\right)=0;$$ hence
it is enough to show that
\begin{equation}\label{goestozero}
\lim_{t\to\infty} \sup_{x\in B} (P\left(\frac{n_t^x}{N e^{\beta
t^*}}<\frac{1}{2} p_t\right)=0.
\end{equation}
Let $R$ denote the law of $N$  and define the conditional laws
$$P^y(\cdot):=P(\cdot\mid N=y),\ y>0.$$
Then
$$P\left(\frac{n_t^x}{N e^{\beta t^*}}<
\frac{1}{2} p_t\right) =\int_0^{\infty}R(\mathrm{d}y)\,
P^y\left(\frac{n_t^x}{y e^{\beta t^*}}<\frac{1}{2} p_t\right).$$
Define the conditional probabilities
$$\widetilde{P}^y(\cdot):= P^y\left(\cdot\mid |Z_{t,x}^B|\ge
\mu_t\right)= P\left(\cdot\mid N=y,\ |Z_{t,x}^B|\ge \mu_t \right),\
y>0,$$ where $\mu_t=\mu_{t,y}:=\lfloor\frac{3y}{4}e^{\beta
t^*}\rfloor$. Recall that (\ref{pt}) defines $p_t$ by taking
supremum over $x$ and  that $|Z^{x,B}_t|$ in fact does not depend on
$x$. One has
\begin{eqnarray}\label{twoterms}
&&P\left(\frac{n_t^x}{N e^{\beta t^*}}< \frac{1}{2} p_t\right)\le \\
&&\int_0^{\infty}R(\mathrm{d}y)\, \left[\widetilde{P}^{y
}\left(\frac{n_t^x}{y e^{\beta t^*}}\nonumber
 <\frac{1}{2} p_t\right)+P^y\left(e^{-\beta t^*}|Z^{x,B}_t|<\frac{3}{4}y\right)\right].
\end{eqnarray}
As far as the second term of the integrand in (\ref{twoterms}) is
concerned, the limit in (\ref{N}) implies that
$$\lim_{t\to\infty}\int_{\mathbb
R}R(\mathrm{d}y)\, P^y\left(e^{-\beta
t^*}|Z^{x,B}_t|<\frac{3}{4}y\right)=
\lim_{t\to\infty}P\left(e^{-\beta
t^*}|Z^{x,B}_t|<\frac{3}{4}N\right)=0.$$

Let us now concentrate on the first term of the integrand in
(\ref{twoterms}). In fact, it is enough to prove that for each fixed
$K>0$,
\begin{equation}\label{withK}
\lim_{t\to\infty} \int_{1/K}^{\infty} R(\mathrm{d}y)\,
\widetilde{P}^{y }\left(\frac{n_t^x}{y e^{\beta t^*}}<\frac{1}{2}
p_t\right)=0.
\end{equation}
Indeed, once we know (\ref{withK}), we can write
\begin{eqnarray}\label{split}
&&\lim_{t\to\infty} \int_{0}^{\infty} R(\mathrm{d}y)\,
\widetilde{P}^{y }\left(\frac{n_t^x}{y e^{\beta t^*}}<\frac{1}{2}
p_t \right)\le \nonumber\\
&&\lim_{t\to\infty} \int_{1/K}^{\infty} R(\mathrm{d}y)\,
\widetilde{P}^{y }\left(\frac{n_t^x}{y e^{\beta t^*}}<\frac{1}{2}
p_t\right)+R\left(\left[0,\frac{1}{K}\right]\right)=
R\left(\left[0,\frac{1}{K}\right]\right).
\end{eqnarray}
Since this is true for all $K>0$, thus letting $K\uparrow\infty$,
\begin{equation*}
\lim_{t\to\infty} \int_{0}^{\infty} R(\mathrm{d}y)\,
\widetilde{P}^{y }\left(\frac{n_t^x}{y e^{\beta t^*}}<\frac{1}{2}
p_t \right)=0.
\end{equation*}
Returning to (\ref{withK}), let us pick randomly $\mu_t$ many
points out of the $J(x,t)$ many particles --- this is almost
surely possible under $\widetilde{P}^{y}$. (Again, `randomly'
means that the way we pick the particles is independent of their
genealogy and their spatial position.) Let us denote the
collection of these $\mu_t$ many particles by $M_t$. Define
$$\widehat{n} _t^x:= \sum_{i\in M_t}1_{A_i},$$ Then, one has
\begin{equation}
 \widetilde{P}^{y}\left(\frac{n_t^x}{y e^{\beta t^*}}<\frac{1}{2} p_t\right)\le
 \widetilde{P}^{y}\left(\frac{\widehat{n}_t^x}{y e^{\beta t^*}}<\frac{1}{2} p_t\right).
\end{equation}
We are going to use Chebysev's inequality and therefore we now
calculate the variance. One has
$$\widetilde{\mathrm{Var}}^y(\widehat{n}_t^x)=
\mu_t(p_t-p_t^2)+\mu_t(\mu_t-1) \frac{\sum_{(i,j)\in
K(t,x)}\mathrm{\widetilde{cov}}(1_{A_{i}},1_{A_{j}})}{\mu_t(\mu_t-1)},$$
where
 $K(t,x):=\{(i,j)\
:\ i\neq j, 1\le i,j\le \mu_t\}.$ Now observe that
$$ \frac{\sum_{i,j\in
K(t,x)}\mathrm{\widetilde{cov}}(1_{A_{i}},1_{A_{j}})}{\mu_t(\mu_t-1)}=\mathbf
{E}\,\mathrm{\widetilde{cov}}(1_{A_{i}},1_{A_{j}})=(\mathbf
{E}\otimes \widetilde{P}^y)(A_i A_j)-p_t^2,$$ where  under
$\mathbf P$ the pair $(i,j)$ is chosen randomly and uniformly over
the $\mu_t(\mu_t-1)$ many possible pairs.

Let $Q^{t,y}$  and $Q^{(t)}$ denote the distribution  of the death
time of the most recent common ancestor of the $i$th and the $j$th
particle under $\widetilde{P}^y$ and under $\widetilde{P}$,
respectively. One has
$$(\mathbf {E}\otimes \widetilde{P}^y)(A_i A_j)=p_t\int_{s=0}^t \int_B
p^t_{t-s,x}\widetilde {p}^{(t)} (0,s,\mathrm{d}x)\,
Q^{t,y}(\mathrm{d}s),$$ where
$$\widetilde {p}^{(t)}(0,t,\mathrm{d}x):=
\mathbb P_0(W_t\in \mathrm{d}x\mid W_z\in B,\ z\le t).$$ By the
Markov property applied at time $s$, $$p_s^t\int_B
p^t_{t-s,x}\widetilde {p}^{(t)} (0,s,\mathrm{d}x)=p_t,$$ and thus
$$(\mathbf {E}\otimes \widetilde{P}^y)(A_i
A_j)=p_t\int_{s=0}^t \frac{p_t}{p_s^t}\, Q^{t,y}(\mathrm{d}s).$$

Hence
\begin{equation}\label{variance}
\widetilde{\mathrm{Var}}^y(\widehat{n}_t^x)\le \mu_t
(p_t-p_t^2)+\mu_t (\mu_t-1)p_t^2\cdot (I_t-1),
\end{equation}
where $$I_t:=\int_{s=0}^{\infty} [p_s^t]^{-1}
Q^{t,y}(\mathrm{d}s).$$ Note that this estimate is uniform in $x$
(see the definition of $p_t$ in (\ref{pt})). Define also
$$J_t:=\int_{s=0}^{\infty} [p_s^t]^{-1}
Q^{(t)}(\mathrm{d}s).$$

 As a next
step, we show that
\begin{equation}\label{lim.of.J}
\lim_{t\to\infty}J_t=1.
\end{equation}
Since $J_t\ge 1$, thus it is enough to prove that
$$\limsup_{t\to\infty}J_t\le 1.
$$
For $r>0$ we denote by $\lambda_r^*:=\lambda_c(\frac{1}{2}\Delta
,B(0,r))$ the principal eigenvalue of $\frac{1}{2}\Delta$ on
$B(0,r)$. Since $\lambda_r^*$  tends to zero as $r\uparrow \infty$
we can pick  an  $R>0$ such that $-\lambda_R^*<\beta$. Let us fix
this $R$ for the rest of the proof.

Let us also fix $t>0$ for a moment. From the probabilistic
representation of the principal eigenvalue (see Chapter 4 in
\cite{P}) we conclude the following: for $\hat\epsilon>0$ fixed
there exists a $T(\hat\epsilon)$ such that for $s\ge
T(\hat\epsilon)$,
\begin{equation*}\log
p_s^t\ge (\lambda_{\rho (t)}-\hat\epsilon)  s.
\end{equation*}
 Hence, for $\hat\epsilon>0$
small enough ($\hat\epsilon<-\lambda_R^*$) and for all $t$
satisfying $\lambda_{\rho(t)}\ge \lambda_R^*+\hat\epsilon$ (recall
that $\lim_{t\to\infty}\rho (t)=\infty$) and $s\ge
T(\hat\epsilon,t)$,
\begin{equation}\label{Lyapunovexponent}\log
p_s^t\ge \lambda_R^* \cdot s.
\end{equation}
Note that $T(\hat\epsilon,t)$ can be chosen uniformly in $t$
because\footnote{In fact $\rho(t)$ can be defined in a way that it
is monotone increasing for large $t$'s.} $\lim_{t\to\infty}\rho
(t)=\infty$, and so we will simply write $T(\hat\epsilon)$.
Furthermore, clearly, $T(\hat\epsilon)$ can be chosen in such a way
that
\begin{equation}\label{T grows}
\lim_{\hat\epsilon\downarrow 0}T(\hat\epsilon)=\infty.
\end{equation}
Depending on $\hat\epsilon$ let us break the integral into two
parts:
$$J_t=\int_{s=0}^{T(\hat\epsilon)} [p_s^t]^{-1} Q^{(t)}(\mathrm{d}s)+
\int_{s=T(\hat\epsilon)}^t [p_s^t]^{-1}
Q^{(t)}(\mathrm{d}s)=:J_t^{(1)} + J_t^{(2)}.$$ We are going to
control the two terms separately.

\medskip\noindent \underline{Controlling $J_t^{(1)}$:}
We show that
\begin{equation}\label{smallerthanone}
\exists\lim_{t\to\infty} J_t^{(1)}\le 1. \end{equation}
 First, it
is easy to check that for all $t>0$, $Q^{(t)}(\mathrm{d}s)$ is
absolutely continuous, i.e.
$Q^{(t)}(\mathrm{d}s)=g^{(t)}(s)\,\mathrm{d}s$ with some
$g^{(t)}\ge 0.$ So
$$J_t^{(1)}=\int_{s=0}^{T(\hat\epsilon)} [p_s^t]^{-1} Q^{(t)}(\mathrm{d}s)
=\int_{s=0}^{T(\hat\epsilon)} [p_s^t]^{-1}
g^{(t)}(s)\mathrm{d}s.$$ Evidently, one has
$[p_{\cdot}^t]^{-1}\downarrow$ as $t\to\infty$. Also, since
$Q^{(t)}([a,b])$ is monotone non-increasing in $t$ for $0\le a\le
b$, therefore $g^{(t)}(\cdot)$ is also monotone non-increasing in
$t$. Hence, by monotone convergence,
$$\lim_{t\to\infty} J_t^{(1)}
=\int_{s=0}^{T(\hat\epsilon)} g(s)\mathrm{d}s=\lim_{t\to\infty}
\int_{s=0}^{T(\hat\epsilon)} g^{(t)}(s)\mathrm{d}s\le
\lim_{t\to\infty} \int_{s=0}^{t} g^{(t)}(s)\mathrm{d}s=1,$$ where
$g:=\lim_{t\to\infty}g^{(t)}$.

\medskip\noindent \underline{Controlling $J_t^{(2)}$:}
Recall that
\begin{equation}\label{Lyap}
\log p_s^t\ge \lambda_R^*\cdot s,\ \forall s\ge T(\hat\epsilon).
\end{equation}
Thus, $$J_t^{(2)}\le \int_{T(\hat\epsilon)}^t
\exp(-\lambda_R^*\cdot s)\, Q^{(t)}(\mathrm{d}s).$$  We will show
that
\begin{equation}\label{doublelimit}
\lim_{\hat\epsilon\downarrow
0}\lim_{t\to\infty}\int_{T(\hat\epsilon)}^t \exp(-\lambda_R^*\cdot
s)\, Q^{(t)}(\mathrm{d}s)=\lim_{\hat\epsilon\downarrow
0}\lim_{t\to\infty}\int_{T(\hat\epsilon)}^t
e^{-(\lambda_R^*+\beta) s}\, e^{\beta s}Q^{(t)}(\mathrm{d}s)=0.
\end{equation}
Recall that $0<\beta+\lambda_R^*$. In order to verify
(\ref{doublelimit}), we will show that given $t_0>0$ there exists
some
 $0<K=K(t_0)$ with the property that
\begin{equation}\label{g.bound}
g^{(t)}(s)\le K s e^{-\beta s},\ \mathrm{for}\ t>t_0,\ s\in
[t_0,t].
\end{equation}
Indeed, it will then follow that
\begin{eqnarray*}
\lim_{\hat\epsilon\downarrow
0}\lim_{t\to\infty}\int_{T(\hat\epsilon)}^t \exp(-\lambda_R^*\cdot
s)\, Q^{(t)}(\mathrm{d}s)=\lim_{T(\hat\epsilon)\to
\infty}\lim_{t\to\infty}\int_{T(\hat\epsilon)}^t
\exp(-\lambda_R^*\cdot s)\, g^{(t)}(s)(\mathrm{d}s)\\\le
K\lim_{T(\hat\epsilon)\to \infty}\int_{T(\hat\epsilon)}^{\infty}
s\,e^{-(\lambda_R^*+\beta) s}\, (\mathrm{d}s)=0.
\end{eqnarray*}
Recall that $Q^{(t)}$ corresponds to the conditional law
$P(\cdot\mid |Z_{t,x}^B|\ge \mu_t)$. We now claim that we can work
with $P(\cdot\mid |Z_{t,x}^B|\ge 2)$ instead of $P(\cdot\mid
|Z_{t,x}^B|\ge \mu_t)$. This is because if $Q_0^{(t)}$ corresponds
to $P(\cdot\mid |Z_{t,x}^B|\ge 2)$, then an easy computation
reveals that for any $\epsilon>0$ there exists a $\hat t_0=\hat
t_0(\epsilon)$ such that for all $t\ge \hat t_0$ and for all $0\le
a<b$,
$$\left|Q^{(t)}([a,b])-Q_0^{(t)}([a,b])\right|\le 2(1+\epsilon) Q_0^{(t)}([a,b]);$$
thus, if
\begin{equation}\label{Q0.bound}
Q_0^{(t)}(\mathrm{d}s)\le L s e^{-\beta s}\mathrm{d}s\
\mathrm{on}\ [t_0,t]\ \mathrm{for}\ t>t_0
\end{equation}
holds with some $L>0$, then also
\begin{equation}\label{Q.bound}
Q^{(t)}(\mathrm{d}s)=g^{(t)}(s)\,\mathrm{d}s\le K s e^{-\beta
s}\mathrm{d}s,\ \ t>t_0\vee \hat t_0,\ s\in [t_0,t]
\end{equation} holds with $K:=L+2(1+\epsilon)$.

The bound (\ref{Q0.bound}) is verified in the appendix.

 \noindent It is now easy to finish the proof of (\ref{lim.of.J}).
 To make the dependence on
$\hat\epsilon$ clear, let us write
$J_t^{(i)}=J_t^{(i)}(\hat\epsilon),\ i=1,2$. Then by
(\ref{smallerthanone}), one has that for {\it all}
$\hat\epsilon>0,$
$$\limsup_{t\to\infty}J_t\le
1+\limsup_{t\to\infty}J_t^{(2)}(\hat\epsilon).$$  Hence,
  (\ref{doublelimit}) yields
$$\limsup_{t\to\infty}J_t\le 1+\lim_{\hat\epsilon\downarrow 0}
\limsup_{t\to\infty} J_t^{(2)}(\hat\epsilon)\le 1,$$ finishing the
proof of (\ref{lim.of.J}).

Once we know (\ref{lim.of.J}), we proceed as follows.  Using
Chebysev's inequality, one has
$$\widetilde{P}^{y}\left(\frac{\widehat{n}_t^x}{y e^{\beta t^*}}<\frac{1}{2}
p_t\right)\le \widetilde{P}^{y}\left(|\widehat{n}_t^x-E^y
\widehat{n}_t^x|>\frac{1}{4} p_t y e^{\beta t^*}\right)\le 16
\frac{\widetilde{\mathrm{Var}}^y(\widehat{n}_t^x)}{p_t^2 y^2
e^{2\beta t^*}}.$$ By (\ref{variance}), we can continue the estimate
by
$$\le 16 \left(\frac{\mu_t (p_t-p_t^2)}{p_t^2 y^2 e^{2\beta t^*}}+
\frac{1}{2}\mu_t (\mu_t-1)\cdot y^{-2}e^{-2\beta t^*}\cdot
(I_t-1)\right).$$ Writing out $\mu_t$, exploiting (\ref{lim.of.J}),
using that the lower limit in the integral is $1/K$, and finally,
dropping the $-p_t^2$ term in the numerator, one obtains that

\begin{equation}\label{vegso}
\int_{1/K}^{\infty}\widetilde{P}^{y}\left(\frac{\widehat{n}_t^x}{y
e^{\beta t^*}}\le \frac{1}{2} p_t\right)\le 12 K p_t^{-1}e^{-\beta
t^*}+ \int_{1/K}^{\infty}R(\mathrm{d}y)\,\frac{1}{2}\mu_t
(\mu_t-1)\cdot y^{-2}e^{-2\beta t^*}\cdot (I_t-1).
\end{equation}
(Recall that $I_t$ in fact depends on $y$.)
  Since $\lim_{t\to\infty}p_t
e^{\beta t^*}=\infty$, thus the first term on the righthand side of
(\ref{vegso}) tends to zero as $t\to\infty$. Recall now that
$\mu_t:=\lfloor \frac{3 y e^{\beta t^*}}{4}\rfloor$. As far as the
second term of (\ref{vegso}) is concerned, it is easy to see that it
also tends to zero as $t\to\infty$, provided one knows
$$\lim_{t\to\infty}\int_0^{\infty} R(\mathrm{d}y)(I_t-1)=0.$$
But $\int_0^{\infty} R(\mathrm{d}y)(I_t-1)=J_t-1$ and so we are
finished by recalling (\ref{lim.of.J}). Hence (\ref{withK})
follows. This completes the proof of the lower estimate in Theorem
\ref{mainthm}. $\hfill\blacksquare$
\section{Some additional problems}
These problems were deferred to the section after the proofs,
because the questions themselves refer to parts of the proof.
\begin{problem}\rm The end of the proof
for the lower estimate in Theorem \ref{mainthm} is basically a
version of the Weak Law of Large Numbers. Using SLLN instead
(and making some appropriate changes elsewhere),  can one get
$$\liminf_{t\to\infty}(\log t)^{2/d}(r_t-\beta)\ge-c(d,\eta)\ \ \ a.s.\ ?$$
\end{problem}
\begin{problem}\rm
The question  investigated in this paper was the (local and global)
growth rate of the population. The next step can be the following:
Once one knows the global population size $|Z_t|$, the model can be
rescaled (normalized) by  $|Z_t|$, giving a population of fixed
weight. In other words, one considers the discrete probability
measure valued process
$$\tilde Z_t(\cdot):= \frac {Z_t(\cdot)}{|Z_t|}.$$
Then the question of the {\it shape} of the population for $Z$
for large times is given by the limiting behavior of the
random probability measures $\tilde Z_t,\ t\ge 0$. (Of course,
not only the particle mass  has to be scaled, but also the
spatial scales are interesting ---  see last paragraph.)

Can one for example locate a {\it unique dominant branch} for almost
every environment, so that the total weight of its complement tends
to (as $t\to\infty$) zero?

The motivation for this question comes from our proof of the
lower estimate for Theorem \ref{mainthm}. It seems conceivable
that for large times the ``bulk'' of the population will live
in a clearing    within distance $\ell(t)$ and with radius
$$\rho(t):=R_0[\log \ell (t)]^{1/d}-[\log\log\ell(t)]^2,\ t\ge
0,$$ where $$\lim_{t\to\infty} \ell(t)=\infty\ \ \mathrm{and}\
\lim_{t\to\infty}\frac{\ell(t)}{t}=0\ \mathrm{but}\
\lim_{t\to\infty}\frac{\log t}{\log \ell(t)}=1.$$

L. Mytnik asked {\it how much the  speed  for free BBM reduces}
due to the presence of the mild obstacle configuration. As we have
seen, for free BBM, the radius of the smallest ball covering the
whole population grows linearly (the velocity is $\sqrt{2\beta}$).
Since $\ell(t)=o(t)$, it would be interesting to know if the speed
in fact becomes sublinear. Note that in the above discussion about
the shape we were only talking about the bulk of the population
and not about {\it individual} particles travelling to very large
distances from the origin.

\end{problem}
\section{Annealed asymptotics of global growth}
In the annealed case we note that the following theorem can be
justified by a method that is similar but much simpler than the
one used to prove Theorem \ref{mainthm}. Let $r_t$ be as in
Theorem \ref{mainthm}.
\begin{theorem} For all $d\ge 1$,
$$
\lim_{t\to\infty}t^{\frac{2}{d+2}}(r_t-\beta)=-\tilde c(d,\nu)\
\mathrm{in}\ (\mathbf{E}\otimes P)-\mathrm{probability}.
$$
(That is, loosely speaking, $r_t\approx \beta-\tilde c(d,\nu)\,
t^{-2/(d+2)}.$)
\end{theorem}
The upper estimate goes exactly the same way as in the proof of
Theorem \ref{mainthm}, by using what we know about the
expectation. The lower estimate however becomes much easier.

Now one does not need the two step `bootstrap' argument as in the
quenched case, and so there is no problem with higher dimensions
(the $d\le 2$ assumption was crucial in the proof for the quenched
case in the first part of the bootstrap method). The only
remaining task is to show that the Law of Large Numbers is in
force for the particle number inside the clearing around the
origin. This can be done essentially the same way as in the last
part of the proof of the quenched case. (The radius of the
clearing is the same as in the classical annealed problem for a
single Brownian particle, cf. the interpretation given for Theorem
\ref{exp.thm}.)
\section{Appendix: Proof of the bound (\ref{Q0.bound})}
We now give the proof of the
bound (\ref{Q0.bound}).  In fact we prove a precise formula for the
distribution of the death time of the most recent common ancestor,
which, we believe, is of independent interest. The result and its
proof are due to W. Angerer and A. Wakolbinger (personal
communication).

For simplicity we set $\beta=1$; the general case is similar. Let
us fix $t>0$. Then for $0<u<t$, one has
\begin{equation}\label{anton}
Q_0^{(t)}(s\le u) = \frac{1 - 2ue^{-u} - e^{-2u} + e^{-t}(2u - 3 +
4e^{-u} - e^{-2u})} {(1 - e^{-t})(1 - e^{-u})^2};
\end{equation}
and so the density is
$$f^{(t)}(u):=\frac{\mathrm{d}Q_0^{t}}{\mathrm{d}\l} (u)=2\,\frac{ e^{-u}(u - 2 + (u + 2)e^{-u}) + e^{-t}(1 -
2ue^{-u} - e^{-2u})} {(1 - e^{-t})(1 - e^{-u})^3},$$ where  $l$
denotes Lebesgue measure on $[0,t]$.

\medskip
\noindent\bf Proof of (\ref{anton}):\rm\ Consider the Yule
population $Y_t:=|Z^B_{t,x}|$ and recall that $Q_0^{(t)}$
corresponds to $P(\cdot\mid Y_t\ge 2)$. The first observation
concerns the Yule {\it genealogy}. Let us pick a pair of
individuals from the Yule population at time $t$, assuming that
$Y_t=j,\ j\ge 2$. Denote by $I$ the size of the population {\it
just before} the coalescence time $s$ of the two ancestral lines
(where `before' refers to backward time). That is, let
$I:=Y_{s+\mathrm{d}t}$. Using some formulae from \cite{EPW05}, we
now show that
\begin{equation}\label{dist.of.I}
P(I = i) = \frac{j + 1}{ j - 1}\cdot\frac{ 2}{ (i - 1)i}\cdot\frac{
i - 1}{ i + 1}.
\end{equation}
 Indeed, the distribution of $I$
equals the conditional distribution of $F$ given ${F \le j}$,
where $F$ is defined as follows. First, the pure birth process
$K=(K_i)$ (with respect to `Yule-time' $i$) is defined in Section
3.5 of \cite{EPW05}, and setting $n=2$,
$$P(K_{i-1}=1\mid K_{i}=2)=\frac{2}{i(i-1)}$$
(see formula $(4.10)$ of the paper). Then the `hitting time' $F$
is defined (in the same section) by $F:=\min\{l:K_l=2\}$. The
formula for the distribution of $F$ (formula $(2.3)$  of the
paper) now becomes
$$P(F\le i)=P(K_i=2)=\frac{i-1}{i+1},\ \hspace{5mm}\ i\ge 2.$$
Let $i\le j$. Then
\begin{eqnarray*}
P(I = i)=P(F = i\mid F \le j)&& = P\big(K_{i-1} = 1,K_i = 2\mid F\le j\big)\\
&&= P\big(K_{i-1} = 1,K_i = 2\mid
K_j=2\big)\\&&=\frac{P\big(K_{i-1} = 1,K_i = 2\big)}{P(K_j=2)}.
\end{eqnarray*}
 From the last three displayed formulae
one arrives immediately at (\ref{dist.of.I}).

Let us now {\it embed}  the `Yule time'  into {\it real time}.
Since a Yule population stemming from $i$ ancestors has a negative
binomial distribution, therefore, using the Markov property at
times $u$ and $u+\mathrm{d}u$, one can decompose
\begin{equation}\label{123}P(Y_u = i - 1, Y_{u+\,\mathrm{d}u} = i, Y_t =
j) =p_1\cdot p_2\cdot p_3,\end{equation} where
\begin{eqnarray*}
&&p_1=e^{-u}(1 - e^{-u})^{i-2},\\ &&p_2=(i-1)\,\mathrm{d}u\ \ \mathrm{and}\\
&&p_3=\binom{j - 1}{ i - 1} e^{-(t-u)i}(1 - e^{-(t-u)})^{j-i}.
\end{eqnarray*}
 That is,
\begin{eqnarray*}&&P(Y_u = i - 1, Y_{u+\,\mathrm{d}u} = i, Y_t =
j) \\&& = (i - 1)\binom{j - 1}{ i - 1} e^{-(t-u)i}(1 -
e^{-(t-u)})^{j-i} e^{-u}(1 - e^{-u})^{i-2}
\,\mathrm{d}u.\end{eqnarray*} Since the pair we have chosen
coalesce independently from the rest of the population, the random
variables $s$ and $I$ are independent. Using the definition of $s$
first and then the independence remarked in the previous sentence,
and finally (\ref{dist.of.I}) and (\ref{123}),
\begin{eqnarray*}
 \ P\big(s\in[u,u+\mathrm{d}u],Y_{s+\mathrm{d}t}=i,Y_t=j\big)
&=&P\big(I = i, Y_u = i - 1, Y_{u+\,\mathrm{d}u} = i, Y_t = j\big)
 \\ &=& P(I = i)P( Y_u = i - 1, Y_{u+\,\mathrm{d}u} = i, Y_t =
 j)\\&=&
 \binom{j - 2} {i - 2} \frac{2(j + 1)} {i(i + 1)}e^{-(t-u)i}\times\\
&& (1 - e^{-(t-u)})^{j-i} e^{-u}(1 - e^{-u})^{i-2}
\,\mathrm{d}u,\end{eqnarray*} for $0<u<t.$

Now, summing from $j = i$ to $\infty$, and from $i = 2$ to
$\infty$, and then dividing by $P(Y_t \ge 2) = 1 - e^{-t},$ one
obtains (after doing some algebra) that for $0<u<t$,
\begin{eqnarray*}Q_0^{(t)}\big(s\in (u, u + \ \mathrm{d}u)\big) =
 \sum_{i=2}^{\infty} e^{-u}\frac {2(2e^{-(t-u)} + i - 1)(1 -
e^{-u})^{i-2}}{ (1 - e^{-t})i(i + 1)}\ \mathrm{d}u&&  \\
=2\cdot\frac{  e^{-u}(u - 2 + (u + 2)e^{-u}) + e^{-t}(1 -
2ue^{-u}-e^{-2u})} {(1 - e^{-t})(1 - e^{-u})^3}\
\mathrm{d}u&.&\end{eqnarray*} Equivalently, in integrated form,
one has (\ref{anton}).\ \ $\hfill\blacksquare$

\bigskip \bf \noindent Acknowledgement.
\rm  I  owe thanks for helpful discussions to the following
colleagues: G. Ben Arous, L. Erd\H{o}s, S. C. Harris, A. E.
Kyprianou, A.-S. Sznitman, B. T\'{o}th and A. Wakolbinger .


\begin{thebibliography}{Poisscat}

\bibitem[AN 2004]{AN} K. B. Athreya and P. E. Ney,  {\em Branching
processes}. Dover, 2004.

\bibitem[AB 2000]{AB} S. Albeverio and L.V. Bogachev. {\em
Branching random walk in a catalytic medium. I. Basic equations.}
Positivity, {\bf 4},   (2000),  41-100.

\bibitem[DF 2002]{DF} D. Dawson and K. Fleischmann {\em Catalytic and mutually
catalytic super-Brownian motions}, in Proceedings
of the Ascona '99 Seminar on Stochastic Analysis,
Random Fields and Applications (R. C. Dalang, M.
Mozzi and F. Russo, eds.), (2002) 89-110.
Birkh\"{a}user, Boston.

\bibitem[DV 75]{DV}  M. Donsker and S.R.S. Varadhan
{\em Asymptotics for the Wiener sausage}  Comm.
Pure Appl. Math., \bf 28\rm\ (1995), 525-565.


\bibitem[E 2000]{E} J. Engl\"{a}nder
{\em On the volume of the supercritical super-Brownian sausage
conditioned on survival,}
Stochastic Process. Appl. \bf 88 \rm (2000), 225--243.

\bibitem[EdH 2002]{EdH} J. Engl\"{a}nder and F. den Hollander {\em  Survival asymptotics for branching Brownian motion in a
Poissonian trap field}, Markov Process.  Related
Fields  \bf 9\rm, No. 3, (2003) 363--389.

\bibitem[EK 2001] {EKpreprint}{ Engl\"ander, J. and Kyprianou, A. E.}
{\em Markov branching diffusions: martingales, Girsanov-type theorems and
applications to the long term behaviour},
Preprint 1206, Department of Mathematics, Utrecht University,
2001, 39 pages. Available electronically at
{\tt http://www.math.uu.nl/publications}


\bibitem[EK 2004] {EK}{ Engl\"ander, J. and Kyprianou, A. E.}
{\em Local extinction versus local exponential
growth for spatial branching processes}, Ann.
Probab. {\bf 32}, No. 1A, (2004)  78--99.

\bibitem[EP99]{EP99} {Engl\"ander, J. and Pinsky, R.} {\em On the construction
and support properties of measure-valued diffusions on
$D\subset R^d$ with spatially dependent branching},  Ann.
Probab. {\bf 27}, No. 2, (1999) 684--730

\bibitem[EPW05]{EPW05}{Etheridge, A., Pfaffelhuber, P. and Wakolbinger, A.}
{\em An approximate sampling formula under genetic hitchhiking},
preprint [ArXiv math.PR/0503485]

\bibitem[KT75]{KT1975}
S.\ Karlin and M.\ Taylor, {\it A First Course in Stochastic
Processes}, Academic Press, New York, 1975.


\bibitem[K 2000] {K} A. Klenke {\em A review on spatial catalytic branching.}  Stochastic
models (Ottawa, ON, 1998),  245--263, CMS Conf. Proc., 26, Amer.
Math. Soc., Providence, RI, 2000.

\bibitem[Ky 2004] {Ky} A. E. Kyprianou
{\em Asymptotic radial speed of the support of supercritical
branching and super-Brownian motion in $R^d$}. To appear in
Markov Process.  Related Fields.

\bibitem[KS 2003] {KS}  H. Kesten, V. Sidoravicius
{\em Branching random walk with catalysts}.
Electron. J. Probab. {\bf 8} (2003), no. 5.
(electronic).

\bibitem[P 1995] {P}  Pinsky, R. G. (1995)
{\em Positive Harmonic Functions and
Diffusion}. \newblock Cambridge University Press.

\bibitem [P 1996] {Pconst}  Pinsky, R. G. (1996) Transience, recurrence and local
extinction properties of the support for supercritical finite
measure-valued diffusions. \textit{Ann. Probab.} 24(1), 237-267.


\bibitem[SLBS 2000] {SLBS1} Shnerb, N. M., Louzoun, Y.,
Bettelheim, E. and Solomon, S. (2000) The
importance of being discrete: Life always wins on
the surface, Proc. Nat. Acad. Sciences {\bf 97} ,
10322-10324.

\bibitem[SLBS 2000] {SLBS2} Shnerb, N. M., Bettelheim, E., Louzoun, Y.,
Agam, O. and Solomon, S. (2001), Adaptation of
autocatalytic fluctuations to diffusive noise,
Phys. Rev. E {\bf 63} , 21103-21108.

\bibitem[Sz 98]{SZ} A. Sznitman
{\em Brownian motion, Obstacles and Random Media.} Springer, 1998.

\end{thebibliography}
\end{document}